# A Finite Dominating Set Approach for the Multi-Item Multi-Period Order Allocation Problem under All-Unit Quantity Discounts and Blending Ratios

Fuhad Ahmed Opu[a], Moddassir Khan Nayeem[a], Hamid Najafzad[a], and Omar Abbaas [a*]

[a] Department of Mechanical, Aerospace, and Industrial Engineering, The University of Texas at San Antonio, San Antonio, TX 78249, USA

**Abstract**

This study addresses the multi-item multi-period order allocation problem under all-unit quantity discounts (AUQD) and blending ratios. A manufacturer makes a single product that requires mixing/assembling multiple ingredients/components with pre-determined blending ratios. We consider a single supplier offering quantity-based discounts which introduces non-linearities to the problem. The objective is to minimize procurement cost which includes purchasing, inventory, and ordering costs. We develop a solution procedure that systematically generates a finite dominating set (FDS) of order quantities guaranteed to include an optimal solution to the problem. A Mixed Integer Linear Programming (MILP) model based on the FDS. Our procedure guarantees optimality and eliminates the need for nonlinear discount modeling. Numerical experiments demonstrate that the proposed MILP achieves optimal solutions with significantly reduced computational effort, up to 99% faster for large-scale instances compared to conventional formulations. Sensitivity analyses reveal that the model dynamically adapts to changes in holding costs, shifting between bulk-purchasing and just-in-time strategies, and identifying cost-sensitive ingredients that drive total system cost.

# 1. Introduction

In today's highly competitive and uncertain business landscape, supply chain management has emerged as a critical strategic function, with purchasing decisions accounting for a significant portion of overall operational costs. Empirical studies show that procurement costs of raw materials and components, in many industries, can comprise 50% to 80% of the final cost of goods sold, highlighting the importance of cost-effective sourcing and order planning [1], [2], [3]. Within this context, Order Allocation (OA) plays a pivotal role, not merely as a byproduct of supplier selection, but as a standalone decision-making problem that directly impacts profitability, inventory dynamics, and supply chain resilience [4]. While supplier selection has long been the focus of purchasing literature, recent studies emphasize that optimal OA, especially in multi-period and complex operational environments, deserves independent attention due to its intrinsic complexity and strategic significance [5].

The complexity of OA increases when a manufacturer produces a product that requires a combination of multiple ingredients or components in its Bill of Materials (BOM). The *blending ratio* in a combination defines the required proportions of each ingredient in the final product. The importance of precise blending ratios is evident in a variety of industries, including food and beverage, pharmaceuticals, chemical manufacturing, and others where product quality and performance depend on following the correct blending ratio. In the coffee industry, for example, blending different coffee beans creates unique flavors. According to a study by the National Coffee Association, 66% of American adults drink coffee daily, with many preferring blended varieties [6]. Leading coffee chains like Starbucks offer more than 30 blends and single-origin premium coffees, and the specific coffee blends (including ratios) are treated as confidential trade secrets [7]. In manufacturing, a product may be an assembly of multiple components. For example, a chair may consist of 1 seat, 1 backrest, 1 base, 5 coasters, 2 handles, and 10-20 screws and bolts. A product missing any of these components or has less than the required number of units of any component is defective. Therefore, ideally inventory levels of the different components at a manufacturing plant should match the required blending ratios. A discrepancy in these levels will cause a stockout of one or more components, creating limiting agents that halt production. The result is inability to deliver finished products while still incurring inventory holding costs for the remaining, unused components.

Suppliers use various strategies to encourage customer orders and align the types and quantities of products in those orders with their operations to achieve maximum efficiency. Among these strategies is discounted pricing. The all-unit quantity discount (AUQD) is one of the most common discount strategies across industries due to its simplicity and greater perceived customer value [5]. Under AUQD, once a specified quantity threshold is met, a reduced price applies to all units purchased [8]. This creates a non-linear pricing structure where the total cost drops at the threshold, incentivizing customers to buy in bulk. Quantity discount mechanisms are common in manufacturing and distribution. For instance, surveys indicate that a majority of purchasing managers (71%) receive AUQD from at least one supplier [9]. In the electronics sector, marketplaces such as Newegg provide quantity-based price breaks where per-unit prices drop once buyers cross specified quantity tiers [10]. AUQD offers benefits for both buyers and suppliers and aligns their interests, improving cost efficiency and strengthening supply chain relationships. It incentivizes larger purchase volumes, enabling suppliers to achieve economies of scale in production, logistics, and administration. This reduces per-unit costs and streamlines operations. Buyers benefit through cost savings on high-volume or frequently purchased items, enhancing profitability and competitiveness.

Beyond financial gains, discount policies also foster long-term supplier-buyer relationships. As highlighted by Pourmohammadreza et al. [5] and Qazi et al. [11], AUQD improves customer satisfaction, reduce costs, and contribute to supply chain resilience. Additionally, such strategies help suppliers differentiate themselves in competitive markets, boosting market share and collaboration opportunities [12], [13].

While AUQD offers numerous advantages, it also introduces complexities that companies must carefully manage. The decision to increase procurement quantities to qualify for discounts must be balanced against inventory holding costs, risks of obsolescence, and production flexibility constraints. For example, in the automotive sector, a manufacturer might receive a 10% discount for ordering more than 10,000 units of a specific tire. However, this bulk purchase creates a risk of surplus inventory if demand falls short of projections, tying up capital and increasing warehousing costs. Similarly, this exposes the firm to obsolescence if design changes occur. In the electronics industry, companies procuring chips must account for the short technological life cycle and demand uncertainty. Holding excess inventory in such cases can lead to costly write-offs and reduced financial flexibility. These challenges are further compounded when the quantity needed to secure a discount does not align with the product's required blending ratio. The real challenge, therefore, lies in aligning order quantities with production needs, inventory considerations, and overall cost-minimization goals.

Multi-period planning environments introduce a new layer of interdependence that changes how firms approach procurement decisions. Unlike single-period OA problems, where firms optimize procurement decisions for an isolated period, multi-period scenarios introduce complex interdependencies between inventory levels, procurement costs, and future order decisions. The cascading effect of past procurement decisions on future cost structures makes myopic OA strategies ineffective, requiring firms to adopt more structured, forward-looking approaches to decision-making. Motivated by these challenges, this study investigates a multi-period procurement problem wherein a manufacturer makes a product that consists of a set of ingredients according to a specific blending ratio. The ingredients are procured from a single supplier who offers AUQD structure. The objective is to determine the optimal order quantities that minimize total procurement, holding, and ordering costs, while satisfying demand requirements.

MILP models and traditional optimization methods offer rigorous exact solutions to procurement problems. However, the combinatorial nature of these problems and the complexities introduced by AUQD and blending ratios can make solving these models computationally expensive [14]. The inherent complexity arises from the exponentially growing decision space associated with multiple ingredients, varying blending ratios, and intertemporal inventory dependencies across procurement periods. To address this issue, we exploit the problem structure to identify a FDS of order quantities guaranteed to include an optimal solution, we call this set the Critical Order Quantities (COQs). The set of COQs is derived from three key sources: supplier discount thresholds, cumulative demand over multiple periods, and residual adjustments accounting for available inventory. Moreover, we develop a MILP formulation that uses the set of COQs to reduce the search space and reach an optimal solution efficiently. Therefore, our proposed approach reduces the solution space without compromising optimality, significantly improving computational tractability.

This study has the following contributions: (1) We formulate a Mixed Integer Nonlinear Programming (MINLP) model that represents the multi-item, multi-period OA problem under AUQD and blending ratio constraints; (2) We identify a FDS of order quantities that is guaranteed to include an optimal solution, thereby reducing computational complexity without compromising optimality; and (3) We develop a computationally efficient COQ-based MILP model and illustrate its effectiveness through numerical experiments, demonstrating its scalability and ability to generate cost-effective procurement plans. Collectively, these contributions address key methodological and computational gaps in the literature and offer practical decision support tools for firms navigating complex sourcing scenarios.

The remainder of this paper is structured as follows. Section 2 provides a comprehensive literature review. Section 3 presents the problem statement. Section 4 details our proposed methodological approach, introducing a traditional MINLP formulation, followed by the derivation of our FDS, and then the COQ-based MILP model. Section 5 presents a series of computational experiments and their results. Finally, Section 6 concludes the paper by summarizing key findings, discussing managerial implications, and highlighting promising avenues for future research.

## 2. Literature Review

Efficient OA in multi-period procurement is complex, especially under the influence of non-linear pricing structures and blending ratios arising from multi-ingredient production requirements. Despite advancements in OA and inventory management, existing literature often treats key complexities in isolation, overlooking their combined impact on procurement efficiency. This section examines the evolution of OA models, the incorporation of quantity discount schemes, and the integration of multi-period inventory management.

### 2.1 OA Models

The OA problem requires several interconnected decisions including selecting which product to procure from which supplier(s), determining the order quantity, and deciding the timing of each order. The decision on "what to purchase and how much" is particularly important as it directly influences the overall cost and supply reliability [15]. This discussion emphasizes two critical dimensions of the problem: whether the procurement involves a single product or multiple products, and whether suppliers offer any form of quantity discount.

*2.1.1 Single-Item, Single-Period Models*

Many foundational contributions in the literature have focused on single-item, single-period OA models, where procurement decisions are made once for a fixed planning horizon. These models are useful for capturing the core dynamics of supplier selection and order distribution without the added complexity of multi-item and multi-period inventory interactions. Depending on the objective, such as minimizing procurement cost, managing supplier risk, incorporating sustainability, or dealing with uncertainty, researchers have employed a wide range of optimization and decision-making tools.

Ghodsypour and O'Brien [16] developed one of the earliest comprehensive frameworks using an MINLP model to allocate a single product across multiple suppliers. Their model integrated various cost components, including unit price, transportation, ordering, and storage, while accounting for supplier capacity and multiple criteria, reflecting the realistic trade-offs procurement professionals often face. Building on the need to evaluate suppliers beyond just cost, Talluri and Narasimhan [17] introduced a model that incorporates variability in supplier performance. Their approach groups vendors based on multiple criteria, such as reliability and delivery performance, allowing buyers to make more informed OA decisions in static, single-period settings.

Faez et al. [18] addressed uncertainty in qualitative supplier evaluation using fuzzy set theory and case-based reasoning. Their mixed-integer programming formulation accounted for capacity constraints and demand satisfaction while modeling the vagueness inherent in decision parameters. This made the approach particularly suited for real-world scenarios where not all supplier attributes are crisply defined. Ruiz-Torres et al. [19] presented a decision-tree-based model that addresses uncertainty by incorporating contingency planning into supplier selection and order allocation decisions. Their approach considers supplier reliability and production flexibility, then develops a base allocation plan alongside state-specific delivery contingencies. This method enhances supply resilience while minimizing fixed, variable, and penalty costs under potential disruption conditions.

Glickman and White [20] include transportation economies, exploring supplier allocation in a distributed retail setting. Their model showed that, in some cases, choosing higher-cost suppliers can be optimal if it enables full truckload shipments, ultimately reducing overall system cost. This work highlights how logistics constraints can shape optimal allocation even when dealing with a single item. Sawik [21] proposed a model focusing on disruption-related risks. The proposed model incorporates supplier protection strategies and emergency inventory pre-positioning. It aims to minimize total costs by optimizing for worst-case disruption scenarios using conditional value-at-risk. From an environmental sustainability standpoint, Shaw et al. [22] addressed carbon emissions as a key factor in supplier evaluation. Using a hybrid Fuzzy Analytic Hierarchy Process (Fuzzy AHP) and fuzzy multi-objective linear programming approach, the model enables procurement professionals to integrate sustainability criteria alongside cost, quality, and delivery performance. This study illustrates how environmental metrics can be effectively embedded into single-period allocation decisions under uncertainty.

Lin [23] introduced a supplier evaluation framework that emphasizes interdependencies among qualitative criteria such as service quality and responsiveness. The proposed multi-criteria decision-making approach provides a structure for selecting the most suitable supplier when the order is made once, and multiple non-cost attributes must be considered simultaneously. Further advancing the integration of sourcing and inventory costs, Mendoza and Ventura [24] developed two single-period MINLP models. One allowed independent order quantities per supplier, while the other required uniform ordering policies. Both models minimized the total procurement cost under supplier capacity and quality constraints, emphasizing the benefit of combining supplier selection with basic inventory control, even in a static setting. To explore decentralized procurement dynamics, Mohammaditabar et al. [25] proposed a game-theoretic analysis of supplier selection and OA. Their study contrasted centralized and decentralized systems using cooperative and non-cooperative frameworks. The result showed that when suppliers act independently, heterogeneous opportunity costs significantly influenced allocation outcomes. In contrast, cooperation can achieve results similar to those in centralized models.

*2.1.2 Multiple-Item, Multiple-Period Models*

While single-item, single-period models have laid the groundwork for OA research, real-world supply chains frequently involve the simultaneous management of multiple items across extended planning horizons. This added dimensionality introduces challenges such as temporal coordination of inventory, fluctuating demand, supplier capacity constraints, and quality variability, all of which complicate the decision of how much to order, for which products, in which periods, and from which suppliers. This subsection reviews key contributions that address these complexities within a multi-item, multi-period OA framework.

Basnet and Leung [26] addressed a deterministic lot-sizing problem involving multiple items, suppliers, and time periods. Their model incorporated product-specific holding costs and supplier-dependent ordering costs, requiring decisions on order quantities, timing, and supplier assignment. An enumerative search algorithm and heuristic were proposed to solve the problem efficiently. Rezaei and Davoodi [27] extended this framework by incorporating imperfect product quality and storage space limitations. Items of substandard quality were discounted and sold before the next delivery. Their model integrated transaction costs, holding costs, and space constraints, and was solved using a Genetic Algorithm. Further, Rezaei and Davoodi [28] developed two multi-objective MINLP models, balancing cost, quality, and service level. One model assumed no shortages, while the other allowed backordering. Notably, they modeled ordering frequency as a cost factor, showing that more flexible policies under backordering can improve cost efficiency.

Osman and Demirli [29] proposed a bi-linear goal programming model for manufacturing aerospace components, focusing on order timing and inventory strategies under demand growth. To solve the complex structure, a modified Benders decomposition method was applied, achieving significant computational gains. Cárdenas-Barrón et al. [30] presented a "reduce and optimize" heuristic for large-scale multi-item, multi-period lot-sizing problems. Their method achieved solutions of comparable quality to those found by CPLEX, but with significantly reduced computational time across benchmark instances, offering a scalable approach to real-world procurement planning.

Addressing demand uncertainty, Kara [31] combined fuzzy decision-making with a two-stage stochastic programming model. While supplier evaluation was included, the core focus was on allocating orders across scenarios, highlighting the value of adaptable procurement strategies under uncertain conditions. Gorji et al. [32] developed a model incorporating imperfect product quality, capital constraints, and supplier capacity. Their MINLP model imposed minimum order quantities and aimed to maximize total profit overtime. A Genetic Algorithm was employed to manage the model's complexity. Türk et al. [33] introduced a two-stage framework integrating supplier risk and inventory planning. In the second stage, a multi-objective evolutionary algorithm was used to allocate orders over multiple periods, minimizing costs while addressing supply risks. They show that their algorithm emerged as the most effective across real and synthetic cases.

## 2.2 Quantity Discount

Quantity discounts (QD) are commonly employed by suppliers as a pricing strategy to incentivize larger orders [34]. In the AUQD scheme, once the purchase quantity exceeds a predefined threshold, the supplier offers a reduced unit price for the entire order, often termed a "target rebate" [35]. Several studies have addressed procurement optimization under AUQD scheme. Kokangul and Susuz [36] proposed a single-item model using AHP and goal programming to balance cost minimization and value maximization under capacity constraints. Razmi and Maghool [37] extended this to a multi-item, multi-supplier context, incorporating capacity and budget constraints with multiple discount types, solving the fuzzy bi-objective problem through the augmented ε-constraint method and heuristics. Similarly, Mazdeh et al. [38] addressed dynamic lot-sizing and supplier selection with and without QD, introducing an extended Fordyce–Webster heuristic to manage the computational complexity of multi-supplier environments.

In the area of integrated shipment and discount planning, Mansini et al. [39] formulated an integer programming model that combines truckload shipping and AUQD policies for multiple items. Lee et al. [40] utilized a Genetic Algorithm with a mixed integer programming model to determine replenishment strategies across multiple suppliers and planning periods under combined all-unit and incremental discounts. Choudhary and Shankar [41] and Ayhan and Kilic [42] focused on multi-objective and integrated decision-making frameworks. Choudhary and Shankar [41] developed a multi-objective linear programming model for supplier selection, lot-sizing, and carrier selection under storage space constraints, employing three variants of goal programming to minimize procurement costs, late deliveries, and rejected items. Ayhan and Kilic [42] combined fuzzy AHP and MILP approaches, where AHP determined the criteria weights and the MILP model optimized supplier selection and OA under AUQD conditions.

Beyond traditional optimization approaches, recent research has explored dynamic procurement mechanisms. Abbaas and Ventura [43], [44] introduced iterative combinatorial auction frameworks for multi-item, multi-sourcing supplier selection and OA problems, where suppliers bid under AUQD schemes. By modeling the problem as MINLPs and applying a buyer's profit-improvement mechanism, their studies demonstrated that fostering competition and dynamic bidding significantly enhances procurement outcomes compared to single-round auctions.

Expanding to broader supply chain perspectives, Tsai et al. [45] integrated AUQD policies into supply chain network design under demand uncertainty, illustrating through mixed-integer models that QD can significantly lower total supply chain costs and improve inventory management. Complementarily, Khan [46] studied a profit optimization inventory problem where demand depends on both selling price and consumption time. Combining AUQD with a prepayment schedule, they derived an optimal inventory and pricing strategy, further showcasing how quantity discount policies can be critical in retailer decision-making under dynamic market conditions.

## 2.3 Research Gaps

Despite significant advancements in OA models, several research gaps persist in current literature. Traditional MILP models struggle with computational effort when applied to multi-period procurement scenarios involving multiple ingredients [30]. The inclusion of complex blending ratios further amplifies the computational burden. As noted by Osman and Demirli [29], solving such problems becomes increasingly intensive, often requiring decomposition methods to maintain tractability. Moreover, while AUQD is widely studied, a limited number of studies integrate these discount schemes with multi-period inventory and procurement planning. Due to the OA complexity, heuristic and metaheuristic approaches are widespread. Few studies demonstrate that their solution methods guarantee optimality under realistic procurement settings. There is a lack of formal frameworks that reduce computational complexity for practically sized problems while preserving solution optimality. Most existing studies either neglect the strategic alignment of discounts with procurement timing or oversimplify inventory holding costs and intertemporal dependencies, thereby failing to provide comprehensive decision support.

To the best of our knowledge, the OA literature rarely considers blending ratios of ingredients to make finished products. Most industries such as food, pharmaceuticals, and auto industry use an array of ingredients mixed or assembled together to make their final products. Existing models do not match consumption and supply rations along with discount schemes and multi-period inventory considerations.

In order to address these limitations this study develops a novel optimization framework. It bridges the gap between computational tractability and guaranteed optimality in complex environments. Our primary contribution is the introduction of COQs, a FDS of economically rational order quantities. We provide theoretical proof that this set is guaranteed to contain an optimal solution. Building on this foundation, we formulate an MILP model. Our approach integrates three key complexities: AUQD, ingredient blending ratios, and multi period inventory dynamics. By leveraging the COQ-based feasible space reduction, our framework reduces the computational limitations of traditional approaches. The result is a scalable and practical model that delivers guaranteed optimal procurement plan.

## 3. Problem Statement

In this study we consider a manufacturing process that involves blending a set of ingredients, denoted by $J$, to make a single finished product. Let $\alpha_j$ be the number of units of ingredient $j$, $j \in J$, needed to produce one unit of the finished product. To fulfill the demand for these ingredients, the manufacturer works with a single supplier, identified through evaluation against predefined economic, environmental, and social criteria. The planning horizon $T$ is defined as a discrete ordered set of time periods. Procurement decisions, including the placement of orders, are made at the commencement of each time period $t$, $t \in T$. Let $q_{j,t}$ be the order quantity of ingredient $j$, at the beginning of time period $t$. The supplier offers an AUQD scheme for the required ingredients. The unit cost of an ingredient varies depending on the order quantity which can qualify for a certain discount level. To define the AUQD scheme, let $N_j$ be the set of discount levels offered by the supplier for ingredient $j$. For each ingredient $j$ and discount level $n$, $n \in N_j$, the supplier specifies a lower bound $l_{j,n}$ for the purchase quantity to qualify for this discount level. It is usually the case that the upper bound of the quantity range of a discount level for an ingredient equals the lower bound of the quantity range of the next discount level for the same item, i.e., $l_{j,n+1}$. The supplier's capacity for an ingredient in $J$, denoted by $u_j$, can serve as the upper bound of the quantity range associated with the highest discount level in $N_j$. The corresponding price per unit at discount level $n$ is denoted by $f_{j,n}$. The price per unit for ingredient $j$ during period $t$ based on the order quantity $q_{j,t}$, is denoted by $c_{j,t}$ and can be expressed as follows:

$$c_{j,t} = \begin{cases} f_{j,1}; & l_{j,1} \leq q_{j,t} < l_{j,2}, \\ f_{j,2}; & l_{j,2} \leq q_{j,t} < l_{j,3}, \\ \ldots \ldots \ldots \ldots \\ f_{j,|N_j|}; & l_{j,|N_j|} \leq q_{j,t} < u_j. \end{cases} \quad (1)$$

This discount structure creates discontinuities in the cost function, as the price per unit drops when a certain quantity threshold is met. The following assumptions are used to formulate and solve the problem. These assumptions are consistent with relevant literature[13], [20], [28], [37].

1. Lead time is constant or negligeable.
2. Shortages are not allowed.
3. Demand is deterministic, known in advance, and happens throughout each period at a constant rate.
4. Unlimited production capacity for the manufacturer. While suppliers share their capacity limit as the upper bound of the highest discount level.
5. The discount schemes offered by the suppliers are identical for the entire planning horizon $T$.
6. The holding cost per unit is constant and does not vary with the level of inventory or item. Additionally, the manufacturer has unlimited storage capacity, meaning that holding costs are only a function of quantity and not influenced by space limitations.

Based on the AUQD scheme and our assumptions we define the purchasing cost as follows:

$$\text{Purchasing Cost } (PC) = \sum_{j \in J} \sum_{t \in T} \sum_{n \in N_j} f_{j,n} \, z_{j,t,n} q_{j,t}, \quad (2)$$

where $z_{j,t,n}$ is a binary variable defined as follows:

$$z_{j,t,n} = \begin{cases} 1, & \text{if discount level } n \text{ is used for ingredient } j \text{ in period } t, \\ 0, & \text{otherwise.} \end{cases} \quad (3)$$

A fixed ordering cost, denoted by $a$, is incurred whenever an order is placed, regardless of the order size or number of ingredients ordered. This ordering cost can be expressed as follows:

$$\text{Ordering Cost } (OC) = \sum_{t \in T} a y_t, \quad (4)$$

where $y_t$ is a binary variable defined as follows:

$$y_t = \begin{cases} 1, & \text{if an order is placed in period } t, \\ 0, & \text{otherwise.} \end{cases} \quad (5)$$

Next, let $d_t$ be the demand of the finished product during time period $t$. Thus, the demand for ingredient $j$ during time period $t$, referred to as $d_{j,t}$, can be written as follows:

$$d_{j,t} = \alpha_j \, d_t. \quad (6)$$

In an ideal case, the proportion between $q_{j,t}$ and $q_{j',t}$ where $j, j' \in J$ should satisfy:

$$\frac{\alpha_j}{\alpha_{j'}} = \frac{q_{j,t}}{q_{j',t}} = \frac{I_{j,t}}{I_{j',t}} = \frac{d_{j,t}}{d_{j',t}}. \tag{7}$$

Satisfying this proportion ensures that the supply of each component matches its consumption. A mismatch in this proportion would lead to a surplus of some ingredients and a shortage of others. This imbalance creates limiting agents that halt production while the company still incurs inventory holding costs for the unused components. However, given the AUQD, the buyer may purchase more than the needed quantity of one or more ingredients to move up to the next discount level which will lower the cost per unit and may offset the additional inventory cost. In this case, the proportions between purchased quantities and inventory levels may not be consistent with Equation (7). In general, let $I_{j,t}$ represent the inventory of ingredient $j$ at the end of time period $t$. $I_{j,t}$ can be defined as follows:

$$I_{j,t} = I_{j,t-1} + q_{j,t} - d_{j,t}. \tag{8}$$

Inventory incurs holding cost, denoted by $h_j$, per unit of ingredient $j$ per time period. Inventory holding cost can be expressed as follows:

$$\text{Holding Cost } (HC) = \sum_{j \in J} \sum_{t \in T} h_j \left( I_{j,t} + \frac{d_{j,t}}{2} \right). \tag{9}$$

Figure 1 shows an example of the inventory level, $I_{j,t}$, for a single ingredient, $j$, over consecutive time periods. At the beginning of each period, if there is an order, the stock level increases upon the arrival of the purchased quantity, $q_{j,t}$. After that, inventory is consumed at a constant rate to satisfy the demand, $d_{j,t}$, as shown by the constant downward slope. The inventory remaining $I_{j,t}$, is carried for the entire period and moves forward to the next period. This carried-over inventory incurs holding costs, representing a potential trade-off where the benefit of a lower purchase price from the AUQD scheme is countered by the cost of holding unused stock.

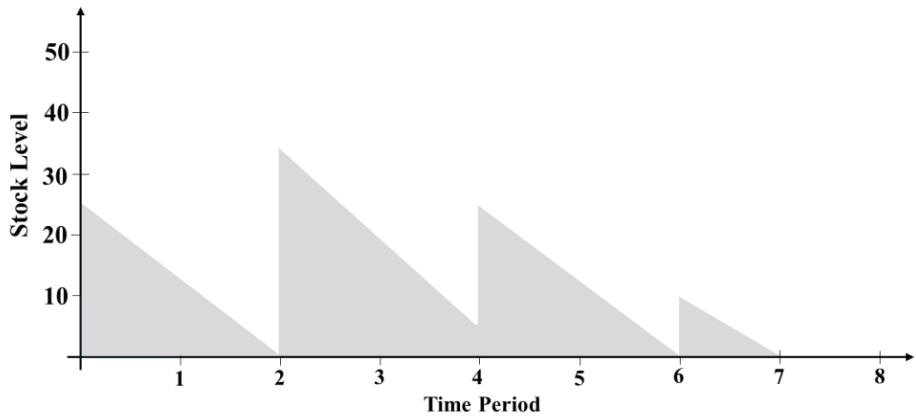

Figure 1: Graphical representation of inventory level for an ingredient over multiple time periods

Figure 2 provides a visual summary of the problem, showing how a single supplier offers distinct AUQD levels for multiple ingredients.

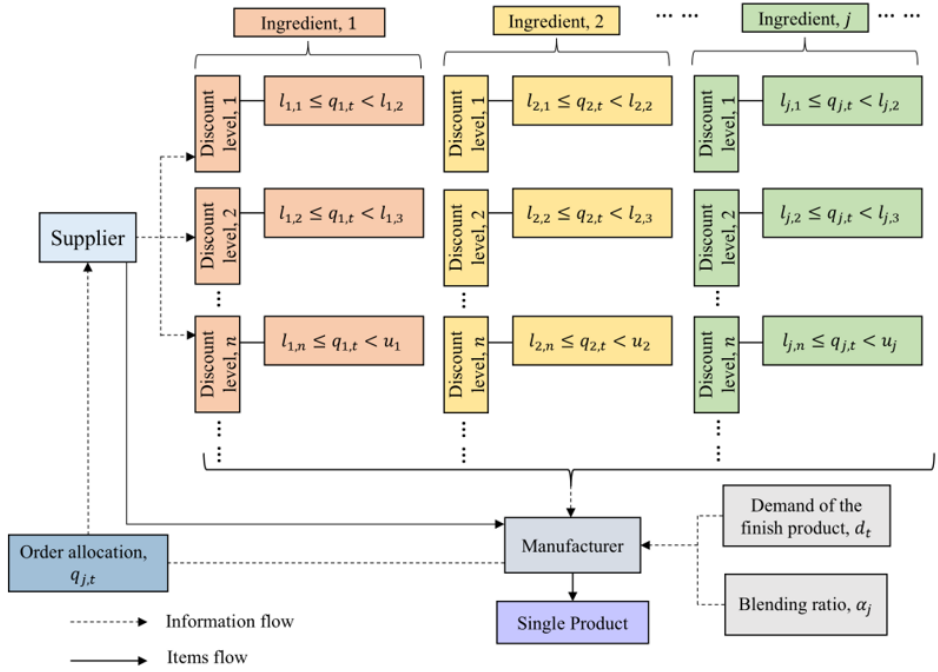

Figure 2: Conceptual framework for multi-period procurement with ingredient blending and nonlinear discount pricing

In this study we will develop an efficient approach to find the optimal order quantity of each ingredient during each time period to minimize the total cost defined as follows:

$$\text{Total Cost } (TC) = PC + OC + HC. \tag{10}$$

## 4. Methodology

In this section, we present our methodological approach to address the procurement optimization problem described in Section 3. We begin by developing a traditional MINLP model serving as our baseline formulation. Recognizing the computational complexity associated with this standard MINLP formulation under multi-period conditions, AUQD, and blending ratios, we introduce the concept of COQs. Following this, we prove that at least one optimal solution exists in this set. Therefore, the set of COQs reduces the search space to an FDS of potential order quantities while ensuring the optimality of the solution. Finally, we integrate these critical order quantities into a MILP model, demonstrating improvements in computational efficiency while preserving solution optimality and practical relevance.

### 4.1 MINLP Model Formulation

Below, we summarize the notation used throughout this model, followed by the mathematical representation of our objective function and constraints.

Sets and indices:
$J$:     Set of ingredients, $j \in J$.
$T$:     Set of time periods, $t \in T$.
$N_j$:    Set of discount levels for ingredient $j$, $n \in N_j$.

Parameters:
$d_t$:     Demand of the finished product in period $t$.
$\alpha_j$:     Number of units of ingredient $j$ needed to produce one unit of the finished product.
$h_j$:     Holding cost per unit of ingredient $j$ per time period.
$f_{j,n}$:   Cost per unit of ingredient $j$ under discount level $n$.
$a$:     Fixed ordering cost per order.
$l_{j,n}$:   Lower bound for discount level $n$ for ingredient $j$.
$u_j$:    Supplier capacity for ingredient $j$.

Decision variables:
$q_{j,t}$:   Order quantity of ingredient $j$ at the beginning of period $t$,
$I_{j,t}$:   State variable determined by the balance equation representing the inventory level of ingredient $j$ at the end of period $t$,
$y_t$:     Binary variable,
$$y_t = \begin{cases} 1, & \text{if an order is placed in period } t, \\ 0, & \text{otherwise.} \end{cases}$$
$z_{j,t,n}$:   Binary variable,
$$z_{j,t,n} = \begin{cases} 1, & \text{if discount level } n \text{ is used for ingredient } j \text{ in period } t, \\ 0, & \text{otherwise.} \end{cases}$$

MINLP model $P$:
$$\text{Minimize: } TC, \tag{11}$$

subject to,

$$I_{j,t} = I_{j,t-1} + q_{j,t} - \alpha_j d_t; \quad \forall j \in J, t \in T, \tag{12}$$

$$l_{j,n} z_{j,t,n} \leq q_{j,t}; \quad \forall j \in J, t \in T, n \in N_j, \tag{13}$$

$$q_{j,t} \leq u_j \sum_{n \in N_j} z_{j,t,n}; \quad \forall j \in J, t \in T, \tag{14}$$

$$\sum_{n \in N_j} z_{j,t,n} \leq y_t; \quad \forall j \in J, t \in T, \tag{15}$$

$$q_{j,t}, I_{j,t} \geq 0; \quad \forall j \in J, t \in T, \tag{16}$$

$$z_{j,t,n}, y_t \in \{0,1\}; \ \forall j \in J, t \in T, n \in N_j. \tag{17}$$

Equation (11) represents the objective function which minimizes the total cost, including purchasing, ordering, and holding costs. These costs are presented in Equations (2), (4) and (9) respectively. Constraint set (12) defines inventory levels at the end of each period which consists of the inventory carried over from the previous period, newly ordered quantity, minus the demand of that period. Constraint set (13) ensures that the order quantity exceeds the lower bound of the selected AUQD level. Constraint set (14) guarantees that order quantities stay within the capacity limit of the supplier. Constraint set (15) ensures that only one discount level is selected and that the manufacturer pays the ordering cost if an order is placed. Constraint set (16) includes non-negativity constraints. Note that the non-negativity of the inventory level variable, $I_{j,t}$, implicitly ensures demand satisfaction in each period. Finally, (17) is the set of domain constraints for binary variables.

### 4.2 Critical Order Quantities (COQs)

Traditional exact optimization methods often face challenges related to computational complexity, especially for practically sized problems. To address this challenge, we derive a FDS of order quantities guaranteed to include an optimal solution to the problem, referred to as COQs. This set reduces the search space for the problem and simplifies the decision-making process without compromising solution optimality. We divide COQs into three categories and introduce them in the following three subsections. After that, we prove that this set is guaranteed to have an optimal solution.

#### 4.2.1 Supplier-Defined Discount Thresholds

The first category of COQs is derived from supplier-defined discount thresholds. Under the AUQD structure, the supplier shares quantity thresholds, $l_{j,n}$, that trigger discounts. Ordering at these thresholds reduces the cost per unit for all purchased units. In a single-item procurement problem, this category will include the lower bounds of the discount levels and the supplier's capacity or the upper bound of the last level as follows,

$$\{l_{j,n} \text{ and } u_j | j \in J, n \in N_j\}. \tag{18}$$

However, in a multi-item problem it is possible that the offered discount levels for the different ingredients do not match the manufacturer's blending ratios. Therefore, the manufacturer may order a quantity that reaches the threshold for a desired discount level for one ingredient and then calculates the required quantity of the other ingredients in $J$ to meet the required blending ratios. Therefore, the set of all discount-based order quantities for ingredient $j$, denoted by $Q_j^T$, is defined as follows:

$$Q_j^T = \{l_{j,n} \text{ and } u_j | n \in N_j\} \cup \left\{ \frac{\alpha_j}{\alpha_{j'}} l_{j',n} \text{ and } \frac{\alpha_j}{\alpha_{j'}} u_{j'} \middle| j' \in J, n \in N_{j'}, \text{and } j' \neq j \right\}. \tag{19}$$

The upper bound of the size of this set can be given using the following formula:

$$|Q_j^T| = |J| + \sum_{j \in J} |N_j|. \tag{20}$$

### 4.2.2 Exact Demand of One or More Consecutive Periods

The second category involves the exact demand of one or more consecutive periods. Considering fixed ordering costs and potential cost benefits from achieving higher quantity discounts, it can be economically advantageous for the manufacturer to consolidate the demands of multiple consecutive periods into fewer, larger orders. We define these aggregated demands as the cumulative sum of demands starting from each period $t$ through all possible future periods $t'$, where $t' \geq t$. We start by calculating this cumulative aggregated demand for the finished product as follows:

$$Q_t^A = \bigcup_{t'=t}^{|T|} \sum_{k=t}^{t'} d_k. \tag{21}$$

Now, given the blending ratios, we can easily find the set of aggregate demand quantities starting from period $t$ for any ingredient $j$ as follows:

$$Q_{j,t}^A = \alpha_j Q_t^A. \tag{22}$$

The set of all possible aggregate demands for ingredient $j$ across the planning horizon is,

$$Q_j^A = \bigcup_{t \in T} Q_{j,t}^A. \tag{23}$$

The upper bound of the number of aggregate demand quantities for one ingredient across all periods can be found using:

$$|Q_j^A| = \sum_{k=1}^{T} (|T| - k + 1) = \frac{|T|(|T|+1)}{2}. \tag{24}$$

### 4.2.3 Residual Adjustments

Under the AUQD scheme, the manufacturer may order quantities higher than the exact demand of one or more consecutive periods in order to qualify for a certain discount level. To facilitate the discussion, let the time period $\hat{t}$ be the period when the order is placed with quantity $q_{j,\hat{t}} \in Q_j^T$. Since quantities in $Q_j^T$ correspond to the supplier imposed discount levels and capacities and not the manufacturer's demand, this can lead to residual inventory carried into a future period, $t \in T$, $I_{j,t-1}$, that is not enough to cover the demand of that period, $I_{j,t-1} < d_{j,t}$. This inventory level can be calculated as follows:

$$I_{j,t-1} = q_{j,\hat{t}} - \alpha_j \sum_{k=\hat{t}}^{t-1} d_k. \tag{25}$$

In this case, an order, with quantity that does not belong to the previous two categories, may be required to complement the residual inventory and fulfill the demand of period $t$,

$$q_{j,t} = d_{j,t} - I_{j,t-1}, \tag{26}$$

or fulfill the demand of period $t$ and a set of consecutive future periods,

$$q_{j,t} = d_{j,t} - I_{j,t-1} + \alpha_j \sum_{k=t+1}^{t'} d_k, \quad t' \in T, t' \geq t + 1. \tag{27}$$

This scenario is illustrated in Figure 3, It shows how a past bulk order, placed to meet a discount threshold, is consumed over time. Eventually, this leads to a residual inventory that is insufficient to meet current demand, which generates the need for a residual adjustment order.

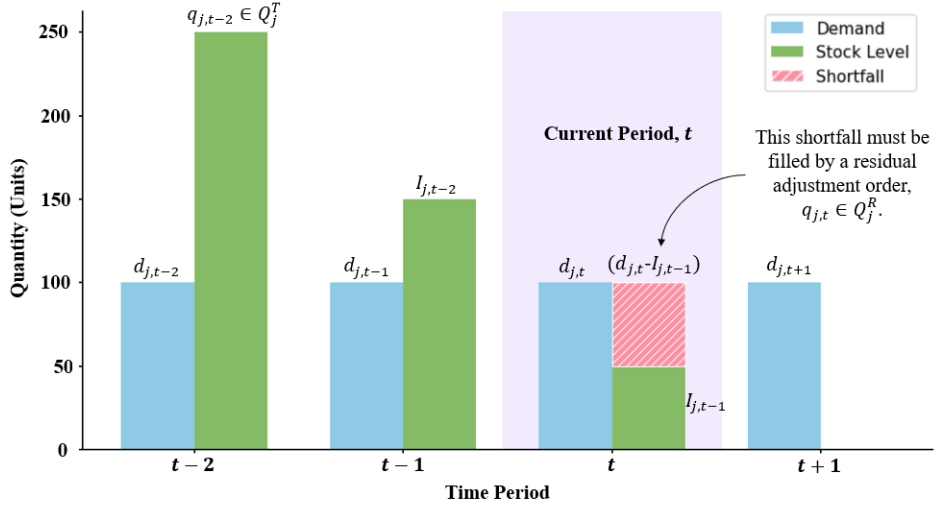

Figure 3: Visual Representation of Residual Adjustment

In general, this case can only happen when the manufacturer is trying to secure a certain discount level by ordering from $Q_j^T$ and then use that quantity to satisfy demand from $Q_j^A$. Let $Q_{j,t}^{A'}$ be a subset of $Q_j^A$ that contains combinations of consecutive demands starting from the first period up to $t$. Mathematically, $Q_t^{A'}$ is defined as follows:

$$Q_t^{A'} = \cup_{\hat{t}=1}^{t} \sum_{k=\hat{t}}^{t} d_k, \tag{28}$$

$$Q_{j,t}^{A'} = \alpha_j Q_t^{A'}. \tag{29}$$

Now, we can find the total set of potential residual inventory adjustments needed to satisfy the demand in period $t$, $d_{j,t}$, using the *Minkowski difference* between $Q_{j,t}^{A'}$ and $Q_j^T$,

$$S_{j,t} = \left\{ b - c \,\middle|\, b \in Q_{j,t}^{A'}, c \in Q_j^T, b - c > 0, b - c < d_{j,t} \right\}, \tag{30}$$

note that $S_{j,t}$ only contains values where $b - c > 0$ and $b - c < d_{j,t}$. Otherwise, there is either no shortage, or the shortage happens in a previous period and needs to be addressed in that period. Now we define $Q_{j,t}^R$ as follows:

$$Q_{j,t}^R = S_{j,t} \cup \{b + c | b \in S_{j,t}, c \in Q_{j,t+1}^A\}. \tag{31}$$

$Q_{j,t}^R$ contains all possible order quantities to complement the inventory carried over and satisfy the demand in period $t$, and potentially some consecutive future periods. Extending this set for ingredient $j$ across the planning horizon we get $Q_j^R$ which is defined as follows:

$$Q_j^R = \bigcup_{t \in T} Q_{j,t}^R. \tag{32}$$

The upper bound of the size of $Q_j^R$ can be found using:

$$|Q_j^R| = \sum_{t=1}^{|T|-1}\left(|Q_j^T|(t-1)(|T|-t)\right) + (|T|-1)|Q_j^T|, \tag{33}$$

which simplifies to:

$$|Q_j^R| = \frac{|Q_j^T|(|T|-1)(|T|^2 - 2|T| + 6)}{6} = |Q_j^T|\left((|T|-1) + \binom{|T|}{3}\right). \tag{34}$$

Next, we combine the three aforementioned categories on the ingredient and period level, and add the possibility of zero order quantity to define the ordered set $Q_{j,t}^{COQ}$, as follows:

$$Q_{j,t}^{COQ} = Q_j^T \cup Q_{j,t}^A \cup Q_{j,t}^R \cup \{0\}, \tag{35}$$

$Q_{j,t}^{COQ}$ is an ordered set where $q_{j,t,k}$ is the $k^{th}$ quantity in $Q_{j,t}^{COQ}$. Finally, we form the overall set of COQs on the finished product level, $Q^{COQ}$, as follows:

$$Q^{COQ} = \bigcup_{t \in T} \bigcup_{j \in J} Q_{j,t}^{COQ}. \tag{36}$$

### 4.2.4 Optimality of COQs

In this subsection, we prove that the set of COQs contains at least one optimal solution making it a FDS.

**Theorem 1.** *There is at least one optimal solution for the multi-item multi-period OA problem under AUQD and blending ratios, where all order quantities belong to $Q^{COQ}$.*

***Proof.*** (By contradiction) Assume there is no optimal solution exists with all order quantities in $Q^{COQ}$. Let $Q^*$ be an optimal solution, and without loss of generality let $q_{j,t}$ be an order quantity that belongs to the optimal solution, $q_{j,t} \in Q^*$, but does not belong to the set of COQs, $q_{j,t} \notin Q^{COQ}, j \in J, t \in T$. By definition, there must be a later period $t' > t$ with $q_{j,t'} > 0$; otherwise $q_{j,t}$ would be equal to the cumulative demand until the end of the planning horizon and it would belong to $Q_j^{COQ}$ leading to a contradiction. Let the AUQD prices per unit for $q_{j,t}$ and $q_{j,t'}$ be $f_{j,n}$ and $f_{j,n'}$, respectively. Also, there must be a small positive number $\epsilon$ such that,

$$\{q_{j,t} \pm \epsilon\} \notin Q^{COQ},$$

where,

$$0 < \epsilon < \min\{|q_{j,t} - \text{nearest value in } Q_{j,t}^{COQ}|, |q_{j,t'} - \text{nearest value in } Q_{j,t'}^{COQ}|\}.$$

Consider the feasible perturbation that moves $\epsilon$ units from period $t$ to $t'$,

$$\tilde{q}_{j,t} = q_{j,t} - \epsilon, \tilde{q}_{j,t'} = q_{j,t'} + \epsilon,$$

with all other components unchanged, the change in $TC$ will be:

$$\Delta TC = \epsilon(f_{j,n'} - f_{j,n}) - \epsilon h_j(t' - t) + 0.$$

Ordering cost does not change since $\epsilon$ is below the smallest change that would create/delete any order. Define the marginal exchange cost $\gamma$ as follows:

$$\gamma = (f_{j,n'} - f_{j,n}) - h_j(t' - t),$$

- If $\gamma < 0$, moving units from $q_{j,t}$ to $q_{j,t'}$ strictly reduces TC, contradicting the optimality of $Q^*$.
- If $\gamma > 0$, moving units in the opposite direction from $q_{j,t'}$ to $q_{j,t}$ strictly reduces TC, contradicting the optimality of $Q^*$.
- If $\gamma = 0$, the cost will not change, one of the two order quantities may move to the nearest COQ to obtain an equally optimal solution with at least one fewer non-COQ quantity.

Repeating this exchange finitely many times yields an optimal solution in which every order quantity lies in $Q^{COQ}$, contradicting our initial assumption and proving that there is at least one optimal solution where all order quantities are in $Q^{COQ}$. □

The theoretical results presented above allow us to use the finite set of COQs without loss of optimality. This approach is expected to reduce the computational complexity, thereby enhancing the practical applicability of our methodology.

### 4.3 COQ-Based Model Reformulation

To incorporate the set of COQ into our mathematical mode, we define the binary variable $x_{j,t,k}$ as:
$$x_{j,t,k} = \begin{cases} 1, \text{if } q_{j,t,k} \text{ is selected,} \\ 0, \text{otherwise.} \end{cases} \tag{37}$$

Correspondingly, $c_{j,t,k}$ represents the purchasing cost associated with $q_{j,t,k}$, given by

$$c_{j,t,k} = \begin{cases} f_{j,1} q_{j,t,k}; & l_{j,1} \leq q_{j,t,k} < l_{j,2}, \\ f_{j,2} q_{j,t,k}; & l_{j,2} \leq q_{j,t,k} < l_{j,3}, \\ \ldots \ldots \ldots \ldots \\ f_{j,|N_j|} q_{j,t,k}; & l_{j,|N_j|} \leq q_{j,t,k} < u_j. \end{cases} \quad (38)$$

The updated purchasing cost formula, $PC'$, can then be expressed as:

$$PC' = \sum_{j \in J} \sum_{k \in K_j} \sum_{t \in T} c_{j,t,k} \, x_{j,t,k}. \quad (39)$$

Replacing Equation (2) by Equation (39) gives the updated total cost formula, $TC'$. Now we introduce our COQ-based MILP as follows:

MILP model $P'$:

$$\text{Minimize}: TC', \quad (40)$$

subject to,

$$I_{j,t} = I_{j,t-1} + q_{j,t} - \alpha_j d_t; \; \forall j \in J, t \in T, \quad (12)$$

$$q_{j,t} = \sum_{k=1}^{|Q_{j,t}^{COQ}|} q_{j,t,k} x_{j,t,k}; \; \forall j \in J, t \in T, \quad (41)$$

$$q_{j,t} \leq y_t M; \; \forall j \in J, t \in T, \quad (42)$$

$$I_{j,t} \geq 0; \; \forall j \in J, t \in T, \quad (43)$$

$$x_{j,t,k}, y_t \in \{0,1\}; \; \forall j \in J, t \in T, k \in [1, |Q_{j,t}^{COQ}|]. \quad (44)$$

Equation (40) represents the objective function which minimizes the total cost. Constraint set (12) defines inventory levels at the end of each period. Constraint set (41) defines the value of $q_{j,t}$ by selecting a member of the set $Q_{j,t}^{COQ}$. Constraint set (42) ensures that the manufacturer pays the ordering cost if an order is placed. Constraint set (43) ensures that stock level is always non-negative. Finally, (44) is the set of domain constraints for binary variables.

## 5 Numerical Experiments

To evaluate the effectiveness and computational performance of the proposed COQ-based MILP model, we conducted a series of numerical experiments. The objectives of these experiments are threefold: (i) to demonstrate the model's ability to generate cost-effective and optimal procurement plans under realistic discount and blending ratio conditions, (ii) to compare the performance of the MINLP formulation with the COQ-enhanced MILP, and (iii) to assess the scalability of the proposed framework across small, medium, and large problem instances.

### 5.1 Test Problem

For the numerical illustration, we consider a test problem involving two ingredients, $j_1$ and $j_2$, required in fixed blending ratios of $\alpha_{j_1} = 3$ and $\alpha_{j_2} = 5$ units, respectively, for each unit of the finished product. The planning horizon spans six periods ($T = 1, \ldots, 6$), with demand sequence (in units of the finished product) $d_t = \{160, 168, 207, 230, 190, 236\}$. Inventory holding costs $h_j$ are set to $1 per unit per period for both ingredients, and the fixed ordering cost is set to $500 incurred whenever an order is placed. Supplier capacity, $u_j$ for both ingredients is assumed to be 5000 units for $j_1$ and 10,000 units for $j_2$, ensuring that supply availability does not constrain the solution. Table 1 summarizes the discount structures for the two ingredients, showing the quantity ranges and corresponding unit costs.

Table 1: Discount thresholds and unit cost

| Ingredients | Discount level, $n$ | Quantity range, $[l_j, u_j)$ | Unit cost ($f_{j,n}$) in $ |
|---|---|---|---|
|  | 1 | $1 \leq q_{j,t} < 1200$ | 15 |
| $j_1$ | 2 | $1200 \leq q_{j,t} < 2500$ | 14 |
|  | 3 | $2500 \leq q_{j,t} < 5000$ | 13 |
|  | 1 | $1 \leq q_{j,t} < 1700$ | 12 |
| $j_2$ | 2 | $1700 \leq q_{j,t} < 4000$ | 10 |
|  | 3 | $4000 \leq q_{j,t} < 10{,}000$ | 8 |

As shown in Table 1, the unit cost for $j_1$ drops from $15 to $14 once the order quantity reaches 1,200 units, and further decreases to $13 when the order reaches 2,500 units. Similarly, $j_2$ has a base price of $12, which falls to $10 at 1,700 units and to $8 at 4,000 units. These discount structures incentivize larger orders but also create the potential for higher inventory holding costs when purchased quantities exceed immediate demand.

## 5.2 Computational Results

Numerical experiments were conducted on a computer with Intel Core i9, 3.0 GHz CPU and 64 GB of RAM. The proposed models were coded using Python 3.11 with the Gurobi Optimizer version 10.0.2.

The COQ-based model identified the optimal procurement strategy for the six-period planning horizon, achieving a total optimal cost of $117,225 and in 0.26 seconds. Table 2 presents the optimal solution. For $j_1$, orders are placed in three periods ($t_1, t_2,$ and $t_3$), while no purchase orders are made in the remaining periods. In the first and second periods, the order quantities cover the exact demand of those periods ($160 \times 3 = 480$, $168 \times 3 = 504$). The order quantity in the third period covers the aggregate demand for the remaining periods (($207 + 230 + 190 + 236) \times 3 = 2589$) and takes advantage of the highest discount level. Therefore, these three order quantities belong to the set $Q_j^A$.

For $j_2$, the solution exploits the AUQD by placing an order of 1,700 units in $t_1$. This order meets the threshold to qualify for the discounted price of $10 per unit, hence $1700 \in Q_j^T$. However, this order leads to an inventory level of 60 units going into period 3 falling short of the demand of $207 \times 5 = 1035$ units. Therefore, an order is required to satisfy the remaining demand of period 3 and potentially the demand of some following periods. The optimal solution orders 4,255 units in $t_3$ to sustains production until $t_6$,

$4255 \in Q_j^R$. As can be seen, optimal order quantities arose from the combination of supplier discount thresholds, cumulative demand aggregates, and residual adjustments which are all COQ. This example validates the theoretical guarantee of optimality within this reduced space.

Table 2: Optimal order quantities and ending inventories

| Ingredients | Metric | $t_1$ | $t_2$ | $t_3$ | $t_4$ | $t_5$ | $t_6$ |
|---|---|---|---|---|---|---|---|
| $j_1$ | $q_{j_1,t}$ | 480 | 504 | 2589 | 0 | 0 | 0 |
| | $I_{j_1,t}$ | 0 | 0 | 1968 | 1278 | 708 | 0 |
| $j_2$ | $q_{j_2,t}$ | 1700 | 0 | 4255 | 0 | 0 | 0 |
| | $I_{j_2,t}$ | 900 | 60 | 3280 | 2130 | 1180 | 0 |

## 5.3 Sensitivity Analysis

To examine the responsiveness of the proposed COQ-based MILP model a comprehensive sensitivity analysis was conducted. This analysis evaluates how the optimal procurement strategy responds to variations in three parameters: (1) inventory holding cost ($h_j$), (2) fixed ordering cost ($a$), and (3) the steepness of the AUQD. For each parameter, "Low" and "High" scenarios were tested against the base case solution presented in Table *2*.

First, the impact of changing inventory holding cost is examined, with results summarized in Table 3. Table 3 presents the optimal procurement plans when the holding cost for both ingredients is varied. The "Low" scenario (where $h_j = \$0.5$) incentivizes order consolidation to reduce ordering cost. This is evident for ingredient $j_1$, where the first two orders from the base case (480 units in $t_1$ and 504 units in $t_2$) are merged into a single, larger order of 984 units in $t_1$. This consolidation saves one fixed ordering cost and is now economically viable due to the low carrying cost.

Table 3: Sensitivity analysis results for inventory holding cost

| Scenario | Ingredients | Metric | $t_1$ | $t_2$ | $t_3$ | $t_4$ | $t_5$ | $t_6$ |
|---|---|---|---|---|---|---|---|---|
| Low | $j_1$ | $q_{j_1,t}$ | 984 | 0 | 2589 | 0 | 0 | 0 |
| | | $I_{j_1,t}$ | 504 | 0 | 1968 | 1278 | 708 | 0 |
| | $j_2$ | $q_{j_2,t}$ | 1700 | 0 | 4255 | 0 | 0 | 0 |
| | | $I_{j_2,t}$ | 900 | 60 | 3280 | 2130 | 1180 | 0 |
| High | $j_1$ | $q_{j_1,t}$ | 480 | 504 | 621 | 1260 | 0 | 708 |
| | | $I_{j_1,t}$ | 0 | 0 | 0 | 570 | 0 | 0 |
| | $j_2$ | $q_{j_2,t}$ | 1700 | 0 | 4000 | 0 | 0 | 255 |
| | | $I_{j_2,t}$ | 900 | 60 | 3025 | 1875 | 925 | 0 |

Conversely, when the inventory holding cost is doubled in the "High" scenario ($h_j = \$2.0$), the solution shifts toward demand-synchronized replenishment to minimize the significant inventory cost. The order frequency for $j_1$ increases from 3 to 5, and the large, consolidated orders from the base case are broken apart into smaller, more frequent procurements. Notably, the solution reduces the order quantity in $t_3$ for $j_2$ from 4,255 units in the base case down to 4,000, the discount threshold, and procures the remaining 255 units in $t_6$. This is done to benefit from the discount while avoiding as much as possible of the high holding cost.

Table 4 shows the change in the optimal solution in response to changing the fixed ordering cost. In the "Low" scenario ($a = 100$), the optimal procurement plan is identical to the base case for both ingredients. This finding suggests that the ordering cost was already a less significant driver than the purchasing and inventory costs, and further reduction was not enough to incentivize a change in strategy.

Table 4: Sensitivity analysis results for fixed ordering cost

| Scenario | Ingredients | Metric | $t_1$ | $t_2$ | $t_3$ | $t_4$ | $t_5$ | $t_6$ |
|---|---|---|---|---|---|---|---|---|
| Low | $j_1$ | $q_{j_1,t}$ | 480 | 504 | 2589 | 0 | 0 | 0 |
| | | $I_{j_1,t}$ | 0 | 0 | 1968 | 1278 | 708 | 0 |
| | $j_2$ | $q_{j_2,t}$ | 1700 | 0 | 4255 | 0 | 0 | 0 |
| | | $I_{j_2,t}$ | 900 | 60 | 3280 | 2130 | 1180 | 0 |
| High | $j_1$ | $q_{j_1,t}$ | 984 | 0 | 2589 | 0 | 0 | 0 |
| | | $I_{j_1,t}$ | 504 | 0 | 1968 | 1278 | 708 | 0 |
| | $j_2$ | $q_{j_2,t}$ | 1700 | 0 | 4255 | 0 | 0 | 0 |
| | | $I_{j_2,t}$ | 900 | 60 | 3280 | 2130 | 1180 | 0 |

However, when the ordering cost is increased in the "High" scenario ($a = 2500$), the solution changes to reduce order frequency. The optimal strategy becomes identical to the low holding cost scenario. For $j_1$, the model consolidates the $t_1$ and $t_2$ orders into a single order of 984 units in $t_1$. This action saves one ordering cost at the expense of a small additional holding cost, confirming the inverse relationship between the two cost parameters.

Finally, the effect of varying the AUQD incentive is shown in Table 5. Here we show the model's sensitivity to the financial incentive offered by the AUQD scheme itself. In the "Low" scenario, the discount steps are minimal. Prices are discounted by $0.5 from one discount level to the next (for $j_1$ $15/$14.5/$14, and for $j_2$ $12/$11.5/$11), compared to the base case of $1 and $2, for $j_1$ and $j_2$ respectively. With this weak incentive, the potential savings from bulk-buying are no longer sufficient to justify the inventory holding costs. The model abandons the base case's consolidation strategy and adopts an approach focused on smaller, more frequent orders. The order frequency for $j_1$ increases from 3 to 4, and for $j_2$ it increases from 2 to 4. The solution stops reaching high discount levels, demonstrating that the marginal savings are not worth the inventory cost. Conversely, in the "High" scenario, where discounts are steep (for $j_1$ $15/$10/$5, and for $j_2$ from $12/$9/$6), the model's behavior is dominated by the pursuit of the lowest possible unit price. For $j_1$, the model places a single large order of 3,573 units in $t_1$ satisfying the demand for the entire planning horizon and forgoing the 3-order strategy of the base case. The savings on unit cost from this single order outweigh the costs of holding inventory for the entire 6-period horizon. Finally, it can be seen that all order quantities in the sensitivity analysis belong to the set of COQ, further demonstrating the validity of our results.

Table 5: Sensitivity Analysis Results for AUQD Discount Steepness

| Scenario | Ingredients | Metric | $t_1$ | $t_2$ | $t_3$ | $t_4$ | $t_5$ | $t_6$ |
|---|---|---|---|---|---|---|---|---|
| Low | $j_1$ | $q_{j_1,t}$ | 480 | 1200 | 0 | 1221 | 0 | 672 |
| | | $I_{j_1,t}$ | 0 | 696 | 75 | 606 | 36 | 0 |
| | $j_2$ | $q_{j_2,t}$ | 800 | 1875 | 0 | 2100 | 0 | 1180 |
| | | $I_{j_2,t}$ | 0 | 1035 | 0 | 950 | 0 | 0 |
| High | $j_1$ | $q_{j_1,t}$ | 3573 | 0 | 0 | 0 | 0 | 0 |
| | | $I_{j_1,t}$ | 3093 | 2589 | 1968 | 1278 | 708 | 0 |
| | $j_2$ | $q_{j_2,t}$ | 1700 | 0 | 4255 | 0 | 0 | 0 |
| | | $I_{j_2,t}$ | 900 | 60 | 3280 | 2130 | 1180 | 0 |

## 5.4 Computational Performance

To further assess the computational efficiency of the proposed COQ-based MILP model ($P'$), we compare its performance with the MINLP model ($P$) across a range of problem sizes. The test instances varied by the number of ingredients, discount levels, and time periods. Table 6 summarizes the solver runtimes (in seconds) for both models obtained using the same computational setup.

Table 6: Comparison of solver runtime for model $P$ and model $P'$ across different problem sizes

| Case | Ingredients | Discount level | Periods | Run Time (sec) | |
|---|---|---|---|---|---|
| | | | | Model $P$ | Model $P'$ |
| 1 | 3 | 3 | 6 | 0.2607 | 0.154 |
| 2 | 3 | 3 | 8 | 3.4105 | 3.3825 |
| 3 | 3 | 3 | 10 | 32.6235 | 6.4861 |
| 4 | 3 | 3 | 12 | 94.7915 | 431.5477 |
| 5 | 3 | 4 | 6 | 0.3256 | 0.2655 |
| 6 | 3 | 4 | 8 | 4.0824 | 3.3511 |
| 7 | 3 | 4 | 10 | 4.7483 | 19.8774 |
| 8 | 3 | 4 | 12 | 854.8181 | 526.399 |
| 9 | 4 | 3 | 6 | 0.3306 | 0.5271 |
| 10 | 4 | 3 | 8 | 8.396 | 4.016 |
| 11 | 4 | 3 | 10 | 230.271 | 13.4503 |
| 12 | 4 | 3 | 12 | 38,229.37 | 179.0712 |
| 13 | 4 | 4 | 6 | 0.4294 | 0.7635 |
| 14 | 4 | 4 | 8 | 9.1501 | 6.2257 |
| 15 | 4 | 4 | 10 | 1,357.763 | 22.9763 |
| 16 | 4 | 4 | 12 | 35,525.38 | 314.9313 |

As shown in Table 6, the COQ-based model ($P'$) achieved faster solution times in 12 out of 16 cases while maintaining global optimality. For small instances, both models were solved almost instantaneously. The performance advantage of our proposed model became more pronounced as problem size increased.

For instance, in case 8 our model was solved in 526.399 seconds, approximately 38% faster than the 854.8181 seconds for *P*. For larger problems, the COQ-based model demonstrated substantial computational savings, solving the largest instance in 314.91 seconds compared to 35,525.38 seconds for model *P*, nearly 99% reduction in runtime. This significant improvement stems from the COQ-based reformulation, which restricts the search space to a finite dominating set of order quantities.

Overall, the COQ-based model not only preserves optimality of the original MINLP formulation but also enhances scalability and computational tractability. These results confirm that the proposed approach can efficiently solve complex, multi-period, multi-ingredient procurement problems that would otherwise be computationally prohibitive using traditional formulations.

## 5. Conclusions

This study addresses a complex multi-period procurement problem involving multiple ingredients that must satisfy fixed blending ratios under an AUQD scheme. In many manufacturing industries, such as food processing, pharmaceuticals, and chemical production procurement managers must align sourcing quantities with blending ratios while navigating nonlinear supplier discounts, inventory dynamics, and multi-period dependencies. These characteristics jointly create a non-convex, computationally intensive optimization landscape that conventional MILP and MINLP models struggle to solve efficiently. To overcome these limitations, this research proposed a COQ-based MILP model that preserves optimality while reducing computational effort.

The concept of COQs represents a FDS of economically rational order quantities derived from three key structural elements: supplier-defined discount thresholds, exact and aggregated period demands, and residual inventory adjustments. Through formal mathematical proof, it was shown that this finite set is guaranteed to contain at least one globally optimal solution to the original MINLP problem. By restricting the decision space to this theoretically validated subset, the proposed COQ-based model effectively eliminates dominated or suboptimal quantity combinations, transforming the nonlinear problem into a tractable MILP without loss of optimality. The results demonstrate that the COQ-based model outperforms traditional formulations in efficiency. For small- and mid-scale instances, the COQ-based model achieved comparable or moderately improved runtimes, with performance gains reaching up to 80% in certain configurations. For larger instances, however, the improvement was substantial, achieving up to 99% reduction of runtimes compared to the conventional MINLP formulation.

The sensitivity analysis reinforced the robustness and interpretability of the COQ-based MILP. By varying holding and ordering cost parameters as well as discount magnitude for different ingredients, the model exhibited consistent and economically logical behavior. For instance, as holding costs increased, total procurement costs rose steadily, accompanied by a strategic shift from large, infrequent orders exploiting quantity discounts to smaller, more frequent orders aligned with near-term demand. Furthermore, higher holding costs led to reduced utilization of higher discount tiers, as maintaining cash flow flexibility and minimizing inventory exposure became more beneficial than securing lower per-unit prices. These findings confirm that the COQ-based MILP not only ensures computational efficiency but also captures realistic managerial trade-offs between inventory holding, discount utilization, and procurement timing.

From a managerial standpoint, the proposed COQ-based MILP framework provides clear and actionable guidance for procurement and production planning within reasonable computational time window. The model demonstrates dynamic policy adaptability, transitioning from bulk purchasing to just-in-time replenishment as cost structure changes. This behavior allows firms to tailor their procurement strategies quickly and accurately without manual policy adjustments. From strategic cost management point of view, firms operating in low-holding-cost environments should exploit supplier discounts through order consolidation, while those facing high storage or obsolescence risks should prioritize smaller, more frequent replenishments. Moreover, the model identifies ingredient-level sensitivities, showing that certain ingredients exert greater influence on total system cost. Managers can therefore target these high-impact inputs for closer monitoring and cost-control efforts. The proposed model guarantees globally optimality and offers reliable decision support in high-stakes industrial settings, ensuring procurement plans that enhance profitability, reduce waste, and strengthen supply chain resilience.

Although the COQ-based MILP developed in this study provides a robust and computationally efficient solution to a challenging class of procurement problems, several promising extensions can further enhance its applicability. First, the current model assumes deterministic demand, which could be relaxed in future work through stochastic or robust optimization techniques to account for demand uncertainty and volatility. A second potential extension involves expanding the framework to a multi-supplier environment with flexible blending ratios or replacement options. In such a setting, procurement decisions would simultaneously consider supplier-specific capacities, reliability indices, and heterogeneous discount structures, while allowing ingredient proportions to vary within feasible bounds based on cost, quality, and availability. This integrated enhancement would capture a more realistic procurement landscape, enabling firms to optimize supplier selection and blending strategies concurrently to manage input volatility and improve cost efficiency.